# K-THEORY. An elementary introduction

## by Max Karoubi

## Conference at the Clay Mathematics Research Academy

The purpose of these notes is to give a feeling of "K-theory", a new interdisciplinary subject within Mathematics. This theory was invented by Alexander Grothendieck[1] [BS] in the 50's in order to solve some difficult problems in Algebraic Geometry (the letter "K" comes from the German word "Klassen", the mother tongue of Grothendieck). This idea of K-theory has invaded other parts of Mathematics, for example Number Theory [Ga], Topology [Bott] (AH1] and Functional Analysis [Connes]. Among many successes of K-theory, one should mention the solution of classical topological problems ([K] chapter V), the Atiyah-Singer index theorem [AS], its wide generalization to the new subject called "Noncommutative geometry" [Connes] and many algebraic applications [MS] [Sr].

This paper is by no means complete and almost no proof is given but the basic ideas are quite simple and no special training is needed to start the subject. However, if one looks for applications in a specific field of Mathematics, some knowledge of this field is required (e.g. Algebraic Geometry, Number Theory, Functional Analysis, etc.).
Our presentation in this paper is driven by one of the most beautiful aspects of the subject which is "Bott periodicity". From its various facets explained here, one guesses some of the main problems and conjectures of the theory.

The prerequisite to read the notes is some familiarity with Topology and elementary Algebra : see for instance the excellent book of Allen Hatcher [Hatcher] or the references below. However, the basic definitions are given in the first section of this paper.

Here are some possible directions for advanced readings. More complete references are listed inside the paper and recalled at the end.

**Introduction to the subject and applications** : [At1], [Ba], [Ber] [Connes], [K1], [M], [R], [S]. One should not be discouraged by unavoidable technicalities and just browse through the main ideas, very often explained at the beginning of these books.
**Some introductory books in Algebraic Topology** : [Hatcher], [ES], [Hu].
**More advanced books and papers** : [FW], [S], [AH1], [AS], [BS], [KV], [L], [Q2], etc.

Finally, and not the least, I would like to thank David Alexandre Ellwood, Research Director of the Clay Mathematics Institure, for helpful comments after a first draft of thie paper.

---

[1] The symbol [ ] means a reference at the end of the paper.

# CONTENTS OF THE PAPER

1. Preliminaries on homotopy theory. Classical Bott periodicity.
2. The beginning of K-theory
3. Relation between K-theory and Bott periodicity
4. K-theory as a homology theory on Banach algebras
5. K-theory as a homology theory on discrete rings

## 1. PRELIMINARIES ON HOMOTOPY THEORY. CLASSICAL BOTT PERIODICITY

Let X and Y be two "nice" topological spaces (for instance metric spaces). Two continuous maps $f_0$ and $f_1$ from Y to X are homotopic if we can find a continuous map[2]

$$F : Y \times [0, 1] \longrightarrow X$$

such that $F(y, i) = f_i(x)$ for $i = 0, 1$. This is an equivalence relation and we denote by [Y, X] the set of equivalence classes. This set is difficult to compute explicitly in general, even for spaces as simple as spheres (see below the partial computation of homotopy groups of spheres).

Two spaces Y and X are homotopically equivalent (or have the same **homotopy type**) if there exist two maps $f : Y \longrightarrow X$ and $g : X \longrightarrow Y$ such that f.g is homotopic to $\text{Id}_X$ and g.f is homotopic to $\text{Id}_Y$. The idea behind this definition is the notion of "shape" of a space which just another word for "homotopy type". Roughly speaking, two spaces have the same shape if they can be deformed one into another. The mathematical meaning of "deformation" is the "homotopy type" which we have just defined.

A variant of these considerations is to take spaces X with a "base point", say $x_0$, and maps preserving base points. For instance, in the definition above of homotopic maps we asssume moreover that $F(y_0, t) = x_0$. One reason for this variant is the definition of the **homotopy groups** of a space X (with a fixed base point). These are defined when Y is a sphere $S^p$ (with a fixed base point) and play an important role in Algebraic Topology.

More precisely, the homotopy classes of base point preserving maps from the sphere[3] $S^p$ to X is the set of path connected components if $p = 0$, form a group if $p > 0$ which is abelian if $p > 1$. This group (or set if $p = 0$) is written $\pi_p(X, x_0)$ or simply $\pi_p(X)$ if the base point is implicit. If two points $x_0$ and $x_1$ are in the same path connected component of X, one can show that $\pi_p(X, x_0)$ and $\pi_p(X, x_1)$ are isomorphic : make a picture ! It is also not difficult to prove that two homotopically equivalent spaces have isomorphic homotopy groups.

---

[2]From now on, we assume that all the maps considered in this paper are continuous.

[3] The sphere $S^p$ is the subset of $\mathbf{R}^{p+1}$ defined by the equation $(x_1)^2 + ... + (x_{p+1})^2 = 1$.

One should add that $\pi_1(X)$ is called the "fundamental group" or "Poincaré group" of a space : Poincaré introduced this group in order to understand the behaviour of the solutions of differential equations around critical points.

The other homotopy groups $\pi_i(X)$ for i > 1 play also a central role in Differential Topology like for example the classification of manifolds. The "**Poincaré conjecture**" for instance is asking the following question. Let M be a topological manifold of dimension 3 which has the homotopy type of $S^3$ ; is M homeomorphic to $S^3$ ? This is one of the problems (with prizes !) raised by the Clay Mathematics Institute.

The homotopy groups of a space are quite difficult to compute. For example, $\pi_p(S^r)$ is not known in general if r > 1. We know however that this group is 0 if p < r and isomorphic to **Z** if r = p. For p > r, some partial information, coming from the fundamental work of Serre, is the following : $\pi_p(S^r)$ is finite if r is odd or if r is even, p ≠ 2r-1. If r is even and p = 2r-1, $\pi_p(S^r)$ is the direct sum of **Z** and a finite group. For these resultas and more general ones, see for instance the basic book of S.-T. Hu on Homotopy Theory [Hu].

However, as it was pointed out by Bott, the homotopy groups of "Lie groups" might be easier to determine in certain circumstances. The Lie groups considered here are defined as the sets of invertible n x n matrices with coefficients in one of the fundamental fields **R**, **C** or **H** : real numbers, complex numbers or quaternions respectively. We denote these groups by $G_n = GL_n(\mathbf{R})$, $GL_n(\mathbf{C})$ or $GL_n(\mathbf{H})$ with the obvious topology (and metric) induced by the natural embedding in a Euclidean space of a suitable dimension.

One can see that $\pi_p(G_n)$ "stabilizes" as n increases. More precisely, the obvious inclusion of $G_n$ in $G_{n+1}$ (adding 0's except for the element at the bottom right corner which is equal to 1) induces an isomorphism

$$\pi_p(G_n) \cong \pi_p(G_{n+1})$$

for sufficiently large n (cf. [K] p. 13-14). More precisely, we have to take n > p+1 (resp. n > p/2, resp. n > p/4 + 1/2) for $GL_n(\mathbf{R})$ (resp. $GL_n(\mathbf{C})$, resp. $GL_n(\mathbf{H})$). We shall denote these stabilized homotopy groups by $\pi_p(GL(\mathbf{R}))$, $\pi_p(GL(\mathbf{C}))$ or $\pi_p(GL(\mathbf{H}))$ respectively. Bott's remarkable result states that :

**THEOREM** [Bott]. *These stabilized homotopy groups are periodic. More precisely*

$\quad\pi_p(GL(\mathbf{R})) \cong \pi_{p+4}(GL(\mathbf{H})) \cong \pi_{p+8}(GL(\mathbf{R}))$ and

$\quad\pi_p(GL(\mathbf{C})) \cong \pi_{p+2}(GL(\mathbf{C}))$

**Remark.** This theorem looks at the first glance quite technical. But, it appeared soon after Bott's discovery that it has a central role in many applications. Topological K-theory (see next sections of the paper) would not have existed without this theorem.

There are other ways to view this theorem and generalize it. If X is a space with base point $x_0$, we call $\Omega X$ the space of base point preserving continuous maps $\sigma$ from the circle $S^1$ to X : this is the "loop space" of X. In the case of a metric space, the topology of $\Omega X$ is induced by the metric $d(\sigma, \tau) = \operatorname{Sup} d_X(\sigma(t), \tau(t))$ where t runs through all the points of the circle. The homotopy groups of $\Omega X$ are those of X with a shift : $\pi_n(\Omega X) \cong \pi_{n+1}(X)$. Moreover, this construction may be iterated : we define for instance $\Omega^2 X$ as $\Omega(\Omega X)$, $\Omega^3(X) = \Omega(\Omega^2(X))$, etc. It is easy to show that $\pi_n(\Omega^r X) \cong \pi_{n+r}(X)$.

With this new definition, Bott periodicity can be stated in a more general form : we have the homotopy equivalences[4]

$$GL(\mathbf{R}) \approx \Omega^4(GL(\mathbf{H})) \approx \Omega^8(GL(\mathbf{R})) \text{ and}$$
$$GL(\mathbf{C}) \approx \Omega^2(GL(\mathbf{C}))$$

If H is a closed subgroup of a Lie group G, the homogeneous space G/H is the quotient of G by the equivalence relation

$$g_1 \sim g_2 \Leftrightarrow (g_2)^{-1} g_1 \in H$$

The "homogeneous space" $GL(\mathbf{C})/GL(\mathbf{R})$ is defined as some kind of limit of the spaces $GL_n(\mathbf{C})/GL_n(\mathbf{R})$ (any real matrix is a complex matrix) ; in the same way $GL(\mathbf{R})/GL(\mathbf{C})$ is the limit of $GL_{2n}(\mathbf{R})/GL_n(\mathbf{C})$, where $GL_n(\mathbf{C})$ is embedded in $GL_{2n}(\mathbf{R})$ via the classical map

$$A + iB \mapsto \begin{pmatrix} A & -B \\ B & A \end{pmatrix}$$

It has also been shown in [Bott] that we have the following homotopy equivalences[5]

$$\Omega(\mathbf{Z} \times BGL(\mathbf{R})) \approx GL(\mathbf{R})$$
$$\Omega(GL(\mathbf{R})) \approx GL(\mathbf{R})/GL(\mathbf{C})$$
$$\Omega(GL(\mathbf{R})/GL(\mathbf{C})) \approx GL(\mathbf{C})/GL(\mathbf{H})$$
$$\Omega(GL(\mathbf{C})/GL(\mathbf{H})) \approx \mathbf{Z} \times BGL(\mathbf{H})$$
$$\Omega(\mathbf{Z} \times BGL(\mathbf{H})) \approx GL(\mathbf{H})$$
$$\Omega(GL(\mathbf{H})) \approx GL(\mathbf{H})/GL(\mathbf{C})$$
$$\Omega(GL(\mathbf{H})/GL(\mathbf{C})) \approx GL(\mathbf{C})/GL(\mathbf{R})$$

---

[4] For a more precise definition of G = $GL(\mathbf{R})$, $GL(\mathbf{C})$ or $GL(\mathbf{H})$, see the end of this section.

[5] See again below a concrete definition of BG when G = $GL(\mathbf{R})$, $GL(\mathbf{C})$ or $GL(\mathbf{H})$.

$$\Omega(GL(\mathbf{C})/GL(\mathbf{R})) \approx \mathbf{Z} \times BGL(\mathbf{R}).$$

[see [K1] p. 148 for a proof and for a concrete definition of the homogeneous spaces G/H in various situations].

For G = GL(**R**), GL(**C**) or GL(**H**), the "classifying space" BG may be explicitly described (up to homotopy) as follows. It is the space consisting of linear continuous maps

$$Q : H \longrightarrow H$$

where H is an infinite Hilbert space (over **R**, **C** or **H**), such that $Q^2 = Q$, the image of Q being a finite dimensional vector space.

In the same way, G has the homotopy type of the space of linear continuous maps

$$\alpha : H \longrightarrow H$$

such that $\alpha - 1$ is a compact operator (i.e. a limit of operators of finite rank).

**Comments and references for further readings.** A nice account of Bott's work can be found in various textbooks and papers, according to one's own taste
1) the original paper of Bott [Bott], however difficult to grasp for a beginner (heavy use of "Morse theory", a deep subject in Riemannian Geometry). We recommend [M2] for a nice overview.
2) an elementary proof (in the complex case) by Atiyah and Bott [AB], also found in the fundamental book of Atiyah [At1]
3) A systematic use of real and complex K-theory (see the following section), where Bott's periodicity appears as a consequence of the periodicity of Clifford algebras [K1].

Concerning Algebraic Topology, there are many good textbooks aimed at undergraduates : we recommend [Hatcher] as the most accessible reference.

**2. THE BEGINNING OF K-THEORY**

The starting point of K-theory was a simple idea of Grothendieck [BS]. If M is a semi-group, one can "**symmetrize**" it in the same way that we define **Z** from the semi-group **N** of natural integers. More precisely, we associate to M a new group S(M) defined as the quotient of M x M by the equivalence relation

$$(a, b) \sim (c, d) \Leftrightarrow \text{there exists e such that}$$
$$a + d + e = b + c + e$$

There is a map from M to S(M) sending m to the class of (m, 0) which is not injective in general. As a matter of fact, the semi-group M is difficult to compute in the applications we have in mind (see the example φ(A) below), whereas S(M) is reachable through standard techniques of Algebra and Topology.

Here is the basic example : let A be a ring with unit and let $M = \phi(A)$ be the set of finitely generated projective (right) A-modules up to isomorphism. Concretely, such a module E is defined as the image of an A-linear map

$$Q : A^n \longrightarrow A^n$$

such that $Q^2 = Q$. If A is a field, E is just a finite dimensional vector space.

As suggested by the notation, M is a semi-group with the following composition law

$$\overline{E} + \overline{F} = \overline{E \oplus F}$$

where $\overline{E}$ denotes in general the isomorphism class of the module E. The group S(M) obtained from M by symmetrization is called the K-theory of A and denoted by K(A).

One should notice that K(A) is a "functor" : a ring homomorphism $f : A \longrightarrow B$ induces a group map $K(A) \longrightarrow K(B)$. On the level of modules, one associates to an A-module E the B-module defined by the formula $F = E \otimes_A B$. If we write E as the image of Q as above, F is the image of the projection operator $\tilde{f}(Q)$ where $\tilde{f}$ is the associated ring homomorphism $M_n(A) \longrightarrow M_n(B)$.

In particular, there is a unique ring map $\mathbf{Z} \to A$. The cokernel of the natural associated map[6] $K(\mathbf{Z}) \longrightarrow K(A)$ is the **reduced K-theory** of A and denoted by $\tilde{K}(A)$. If A is commutative, K(A) is isomorphic (non naturally) to the direct sum $\mathbf{Z} \oplus \tilde{K}(A)$. Therefore, we may view $\tilde{K}(A)$ as the non trivial part of the K-theory of the ring A.

Despite the apparent simplicity of this definition, K-theory has been very useful in many parts of Mathematics. Two extreme examples will show the versatility of the notion.

The first example comes from Number Theory : A is the quotient ring $\mathbf{Z}[x]/(1 + x + .. + x^{p-1})$, where p is an odd prime number. Then $\tilde{K}(A)$ is a finite group studied in the 19[th] century by Kummer (under a different name) and related to Fermat's last theorem : if this group has no p torsion, Kummer proved that Fermat's equation $X^p + Y^p = Z^p$ has no integral solution (except if X or Y = 0 of course). See for instance [BorS] and [M1] where reduced K-theory is related to **ideals in number fields** (i.e. finite extensions of **Q**).

The second example comes from Topology : A is the ring of (real or complex) continuous functions on a compact space X. We shall denote this example by $C_\mathbf{R}(X)$ or $C_\mathbf{C}(X)$ in the real or complex case respectively. Then we have an isomorphism between $K(A) \otimes_\mathbf{Z} \mathbf{Q}$ and some part of the **cohomology**[7] of X (see [K3] for a proof). More precisely :

---

[6] It is relatively easy to show that $K(\mathbf{Z})$ is isomorphic to **Z** ; see the table of K-groups below.

$$K(C_{\mathbf{C}}(X)) \otimes_{\mathbf{Z}} \mathbf{Q} \cong \bigoplus_k H^{2k}(X; \mathbf{Q})$$

$$K(C_{\mathbf{R}}(X)) \otimes_{\mathbf{Z}} \mathbf{Q} \cong \bigoplus_k H^{4k}(X; \mathbf{Q})$$

One should also notice that the K-theory of these rings is a "**homotopy invariant**" : replacing X by X x [0, 1] does not change the K-theory.

The following table lists some basic examples (with comments) one can try to compute. Calculating "K-groups" is the best way to learn about the "theory" behind K-theory.

**Examples of computations (with comments)**

| A | K(A) |
|---|---|
| **Z** | **Z** |
| | This means that every abelian group which is a direct summand of some $\mathbf{Z}^n$ is itself of the type $\mathbf{Z}^m$ for a certain m. |
| $\mathbf{Z}(\sqrt{-5})$ | $\mathbf{Z}/2 \oplus \mathbf{Z}$ |
| | The factor **Z**/2 is related to the existence of a non principal ideal [Sa]. |
| **C**[G] | $\mathbf{Z}^N$ |
| | Here G is a finite group, N the number of its irreducible representations. This computation is the starting point of the beautiful theory of representations of finite groups [Se]. |
| $C_{\mathbf{C}}(S^1)$ | **Z** |
| | The geometric interpretation of this result is the following : any complex vector bundle over the circle is trivial ([K] p. 13) |
| $C_{\mathbf{R}}(S^1)$ | $\mathbf{Z}/2 \oplus \mathbf{Z}$ |
| | The factor **Z**/2 is due to the existence of the Mœbius band ([K1] p. 13). |
| $C_{\mathbf{C}}(S^2)$ | $\mathbf{Z} \oplus \mathbf{Z}$ |
| | One of the factor **Z** is due to the existence of a "canonical" complex line bundle L over the sphere of dimension 2 ([K1] p. 13 or [K2] p. 22). |
| $C_{\mathbf{R}}(S^2)$ | $\mathbf{Z}/2 \oplus \mathbf{Z}$ |
| | The underlying real bundle to L - let us call it L' - is not trivial and is responsible for the factor **Z**/2 since L' $\oplus$ L' is a trivial real bundle ( use [K1] p. 13 again). |

---

[7] The "cohomology" is defined for instance in [Hatcher] or [ES].

| | |
|---|---|
| $M_n(A)$ | $K(A)$ |

This means that K-theory does not change when A is replaced by the ring of n x n matrices with coefficients in A (Morita invariance). Many elementary proofs are possible ; see for instance [K1] p. 138.

| | |
|---|---|
| $\mathcal{K}$ | $\mathbf{Z}$ |

Here $\mathcal{K}$ denotes the algebra of compact operators in an infinite dimensional Hilbert space H. The proof requires some analysis. See for instance [K] p. 176 (exercise 7.9).

| | |
|---|---|
| $\mathcal{B}(H)$ | $0$ |

$\mathcal{B}(H)$ denotes the algebra of all bounded operators in H. The proof is quite elementary as soon we notice that the category of finitely generated projective modules over this algebra is equivalent to the category which objects are Hilbert spaces (direct summands of some $H^n$) and morphisms linear continuous maps between them.

## Comments and recommendations for further readings

There are many books on the elementary aspects of K-theory. The interested student should read for instance [At1] or the first two chapters of [K1] (try to solve the relevant exercises p. 44-51 and p. 105-111).

## 3. RELATION BETWEEN K-THEORY AND BOTT PERIODICITY

We shall limit ourselves to the complex case for simplicity and describe an explicit map from $\pi_{p-1}(GL(\mathbf{C}))$ to $\tilde{K}(A)$, where A is $C_\mathbf{C}(S^p)$, the ring of complex continuous functions defined on the sphere $S^p$. For this purpose, we decompose the sphere $S^p$ into two "fat" hemispheres $S^p_+$ and $S^p_-$, each hemisphere being defined by $x_{p+1} > -1/2$ (resp $x_{p+1} < 1/2$).

Using the last coordinate $x_{p+1}$ as a parameter, we define two continous functions $\alpha_+$ and $\alpha_-$ from the sphere $S^p$ to the unit interval [0, 1] such that $\alpha_+$ (resp. $\alpha_-$ ) is zero outside $S^p_+$ (resp. $S^p_-$) and such[8] that $\alpha_+ + \alpha_- = 1$. We define

---

[8] such a data is called a partition of unity associated to the covering of the sphere by its two hemispheres (see [Mu] p. 222 in a much broader context).

$$\beta = \frac{\alpha}{\sqrt{(\alpha_+)^2 + (\alpha_-)^2}}$$

where $\alpha = \alpha_+$ or $\alpha_-$ and $\beta = \beta_+$ or $\beta_-$ accordingly.

After these preliminaries let us consider an element of $\pi_{p-1}(GL(\mathbf{C}))$ represented by a continuous map (for large enough n)

$$f : S^{p-1} \longrightarrow GL_n(\mathbf{C})$$

We extend this map to $S^p_+ \cap S^p_-$ by meridian projection on the equator

$$\tilde{f} : S^p_+ \cap S^p_- \longrightarrow GL_n(\mathbf{C})$$

The matrix

$$Q = \begin{pmatrix} (\beta_+)^2 & \beta_+\beta_- \, \tilde{f} \\ \beta_+\beta_- \, \tilde{f}^{-1} & (\beta_-)^2 \end{pmatrix}$$

(where $(\beta_+)^2$, $(\beta_-)^2$ or $(\beta_+\beta_-)$ means this scalar times the identity matrix of order n) is a 2n x 2n matrix whose entries are continuous functions on $S^p$ with values in $M_{2n}(\mathbf{C})$. It is easy to see that $Q^2 = Q$. Therefore, the image of Q defines a projective module over the ring $A = C_{\mathbf{C}}(S^p)$, hence an element of $\tilde{K}(A)$ which is independent of the choice of the partition of unity $(\beta_+, \beta_-)$ ( this can be proved by a homotopy argument).

**THEOREM.** *The correspondence* $f \mapsto \text{Im}(Q)$ *defines an isomorphism*

$$\pi_{p-1}(GL(\mathbf{C})) \xrightarrow{\approx} \tilde{K}(C_{\mathbf{C}}(S^p))$$

[for a proof, see for instance [K1] p. 36 and 55 §1.18].

The previous theorem gives us a translation of Bott periodicity in terms of K-theory : one simply has to show that $\tilde{K}(C_{\mathbf{C}}(S^p))$ is periodic of period 2 with respect to p. Equivalently, using the same construction in the real case, one has to show that $\tilde{K}(C_{\mathbf{R}}(S^p))$ is periodic with respect to p of period 8.

**4. K-THEORY AS A HOMOLOGY THEORY ON BANACH ALGEBRAS**

So far, we have defined K-theory for unital rings only. A slight modification allows us to generalize this definition to **non unital rings** A as follows. We assume that A is a k-algebra where k is any commutative ring with unit (for instance **Z**). We define a new <u>unital</u> ring $A^+$ as A x k with the obvious addition and the following "twisted" multiplication

$$(a, \lambda)(a', \lambda') = (aa' + \lambda a' + \lambda' a, \lambda \lambda')$$

There is an obvious "augmentation" $A^+ \longrightarrow k$ and the K-theory of A is then defined as the kernel of the induced map

$$K(A^+) \longrightarrow K(k)$$

It can be shown (not quite easily) that this definition is in fact independent of k.

An interesting and motivating example is the ring of k-valued continuous functions f on a locally compact space X (k = **R** or **C**) such that f(x) goes to 0 when x goes to ∞. Then $A^+$ is the ring of continuous functions on the "one point compactification" of X. For instance, if X = $\mathbf{R}^p$, K(A) is isomorphic to $\tilde{K}(C(S^p))$.

Also note that this method enables us to define a morphism $K(A) \longrightarrow K(B)$ each time we have a general ring map $A \longrightarrow B$ (we no longer assume that f(1) = 1, even if A and B are unital rings).

A sequence of rings and maps

$$0 \longrightarrow A' \longrightarrow A \longrightarrow A'' \longrightarrow 0$$

is called exact if the underlying sequence of abelian groups is exact. In other words A' is a two sided ideal in A, whereas A" may be identified with the quotient ring A/A'.

**THEOREM**. *The previous sequence induces an exact sequence of K-groups*

$$K(A') \xrightarrow{\alpha} K(A) \xrightarrow{\beta} K(A'')$$

*that is* Im($\alpha$) = Ker($\beta$) [for a proof see for instance [M] or [KV]).

It is natural to ask what is Ker($\alpha$) and Coker($\beta$). If one is familiar with homological algebra, one should define "**derived functors**" $K_n$, n ∈ **Z**, of the K-group in order to extend the previous exact sequence to the left and to the right (the group $K_0(A)$ being K(A)). A partial solution (to the left) is given by the theorem a few lines below.

There are at least two ways to solve this problem of derived functors. The first one is to put some topology on the rings involved (i.e. consider Banach algebras as explained below), The other is to stay in pure algebra, which is paradoxically much harder. We shall begin with the first approach.

We recall that a **Banach algebra** (over k = **R** or **C**) is a k-algebra A (not necessarily unital) with the following properties :

1) A norm a $\mapsto$ $\|a\|$ is defined on the vector space A in such a way that A is complete for the distance d(a, b) = $\|a-b\|$
2) We have the inequality $\|ab\| \leq \|a\|.\|b\|$

A typical example is the ring of continuous functions f on a locally compact space X such that f(x) goes to 0 when x goes to $\infty$. The norm $\|f\|$ of f is then the maximum of the values of $|f(x)|$ when x $\in$ X.
More generally, if A is a Banach algebra and X a locally compact space, we define a new Banach algebra A(X) as the ring of continuous functions f on X with values in A with the same condition at infinity. For the definition of the norm, we just replace $|f(x)|$ in the previous example by $\|f(x)\|$. With obvious notations, we have the following isomorphism A(X)(Y) $\cong$ A(X x Y).

**THEOREM.** *Up to isomorphism, there is a unique way to define functors* $K_n(A)$, *for a Banach algebra* A *and* n $\geq$ 0, *such that the following axioms are satisfied :*

**1) Exactness**. *Given an exact sequence of Banach algebras*

$$0 \longrightarrow A' \longrightarrow A \longrightarrow A'' \longrightarrow 0$$

*such that* A' *has the induced norm and* A" *the quotient norm, we have a long exact sequence (infinite to the left)*

$$\longrightarrow K_{n+1}(A'') \longrightarrow K_n(A') \longrightarrow K_n(A) \longrightarrow K_n(A'') \longrightarrow K_{n-1}(A') \longrightarrow$$

**2) Homotopy invariance**. *The canonical inclusion of* A *into* B = A([0, 1]) *induces an isomorphism*

$$K_n(A) \cong K_n(B)$$

**3) Normalization**

$$K_0(A) = K(A)$$

[for a reference, see for instance [K1] p. 109, exercise 6.14].

It is not too difficult to deduce from these axioms the following isomorphism

$$K_n(A) \cong K(A(\mathbf{R}^n))$$

In particular, $\tilde{K}(C_\mathbf{R}(S^n)) \cong K_n(\mathbf{R})$ and $\tilde{K}(C_\mathbf{C}(S^n)) \cong K_n(\mathbf{C})$. Moreover, if A has a unit, the group $K_n(A)$ may also be defined as a suitable homotopy group : for $n > 0$, this is just the homotopy group $\pi_{n-1}(GL(A))$ which is a "limit" of the homotopy groups $\pi_{n-1}(GL_r(A))$ when r goes to ∞.

The previous theorem enables us to define a "cohomology theory" on the category of compact spaces in the sense of Eilenberg and Steenrod [ES] : if X is a compact space and Y a closed subset, we define $K^{-n}(X, Y)$ as the $K_n$ group[9] of the Banach algebra $C(X - Y)$. There are in fact two theories involved according to whether we consider real or complex continuous functions on X - Y. These groups satisfy all the axioms of a cohomology theory [ES], except for the dimension axiom. This cohomology theory was defined initially by Atiyah and Hirzebruch [AH1] and was the starting point of "**topological K-theory**" for compact spaces and later on for Banach algebras [K]. This theory for Banach algebras (developped in the last 20 years) got an enormous importance in a new field of Mathematics called "Noncommutative Geometry" [Connes]. This shows the importance of the following theorem in this new field :

**THEOREM** ([K1] § III) . *Let* A *be a complex (resp. real) Banach algebra. Then we have an isomorphism*
$$K_n(A) \cong K_{n+2}(A) \quad (\text{resp. } K_n(A) \cong K_{n+8}(A))$$

**Remark.** If we stay in the complex case for simplicity and if $n \in \mathbf{Z}$, we see that we can define formally $K_n(A)$ as $K_{n+2r}(A)$ for r sufficiently large such that $n + 2r > 0$. In other words, we can define $K_{-1}(A)$ as $K_1(A)$, $K_{-2}(A)$ as $K_0(A)$, etc. Therefore, to an exact sequence of Banach algebras as above

$$0 \longrightarrow A' \longrightarrow A \longrightarrow A'' \longrightarrow 0$$

we can associate an "hexagonal" exact sequence

$$\begin{array}{ccccc} K_0(A') & \longrightarrow & K_0(A) & \longrightarrow & K_0(A'') \\ \uparrow & & & & \downarrow \\ K_1(A') & \longrightarrow & K_1(A) & \longrightarrow & K_1(A'') \end{array}$$

This solves the problem we had before about the kernel of the map $K_0(A') \longrightarrow K_0(A)$ AND the cokernel of the map $K_0(A) \longrightarrow K_0(A'')$.

There are many proofs of Bott periodicity in the literature. They all require some sophistication. One conceptual proof is to introduce "cup-products" in K-theory. More

---
[9] Note that $K^n(X, Y)$ has not yet been defined for $n > 0$.

precisely, a bilinear pairing between three unital rings A, B and C is given by a **Z**-bilinear map

$$\varphi : A \times C \longrightarrow B$$

such that

$\varphi(aa', cc') = \varphi(a, c) \varphi(a', c')$

$\varphi(1, 1) = 1$

It is easy to see that such a pairing induces a group homomorphism

$$K(A) \otimes K(C) \longrightarrow K(A \otimes_\mathbf{Z} C) \longrightarrow K(B)$$

This definition may be generalized easily to non unital rings. Moreover, if we assume that A, B and C are Banach algebras and that $\|\varphi(a, c)\| \leq \|a\| \cdot \|c\|$, such a bilinear map $\varphi$ induces another one

$$A(\mathbf{R}^n) \times C(\mathbf{R}^p) \longrightarrow B(\mathbf{R}^{n+p})$$

and therefore a bilinear pairing (called the "cup-product")

$$K_n(A) \times K_p(C) \longrightarrow K_{n+p}(B)$$

Another version of Bott periodicity is then the following :

**THEOREM** (see for instance [K4] § 2)
a) *Let* $C = \mathbf{C}$ *and let* $A = B$ *be a complex Banach algebra. Then the group* $K_2(\mathbf{C})$ *is isomorphic to* **Z** *and the cup-product with a generator* u *induces an isomorphism*

$$\beta : K_n(A) \longrightarrow K_{n+2}(A)$$

b) *Let* $C = \mathbf{R}$ *and let* $A = B$ *be a real Banach algebra. Then the group* $K_8(\mathbf{R})$ *is isomorphic to* **Z** *and the cup-product with a generator induces an isomorphism*

$$K_n(A) \longrightarrow K_{n+8}(A)$$

## Sketch of a proof of the theorem, part a)

It relies on the construction of a full theory $K_n(A)$ defined[10] for all $n \in \mathbf{Z}$ and satisfying the formal properties of a cup-product (for instance associativity). We then construct an element v in $K_{-2}(\mathbf{C})$ such that $u \, v = v \, u = 1$ in $K_0(\mathbf{C}) = \mathbf{Z}$. This element v is used in order to construct a morphism "going backwards"

---

[10] In order to avoid a vicious circle, we DONT assume periodicity to define $K_n$ for $n < 0$.

$$\beta' : K_{n+2}(A) \longrightarrow K_n(A)$$

It is simply defined by $\beta'(x) = x.v$. The fact that $\beta$ and $\beta'$ are inverse to each other is a formal exercise. For instance, $(\beta\beta')(x) = \beta(x.v) = (x.v).u = x.(v.u) = x.1 = x$.

**Applications**. A variety of applications of Bott periodicity and K-theory may be found in various papers and books [At1], [AH2], [AS2] [Connes]. We just mention two topological applications, due to J.F. Adams (the proofs may be found in [K], chapter V)

1) The sphere $S^n$ is a topological group only if $n = 0, 1$ or $3$
2) What is the maximum number $\rho(t)$ of linearly independent tangent fields on the sphere $S^{t-1}$ ?
    a) According to an old theorem of Hopf, $\rho(t) = 0$ if t is odd
    b) If t is even, we can write it as (odd number).$2^\beta$ and we divide $\beta$ by 4, say $\beta = \gamma + 4\delta$ where $0 \leq \gamma \leq 3$.
    Then $\rho(t) = 2^\gamma + 8\delta - 1$

**Suggestions for advanced readings** : the books [At1] and [K], already quoted many times.

## 5. K-THEORY AS A HOMOLOGY THEORY ON DISCRETE RINGS

For various reasons, especially the applications of K-theory to Algebraic Geometry and Number Theory, the definition of the functors $K_n$ we have just given is not very satisfactory. We would rather not use the topology of the ring A. For such a purpose, a definition has been proposed by Quillen in 1970 [Q2]. Unfortunately, Quillen's definition requires some sophistication in Algebraic Topology. Therefore, we shall use another one, introduced slightly before by Villamayor and the author [KV] which coincides with Quillen's definition in many favourable cases. It is much easier to define, in the same spirit as the higher K-groups for Banach algebras.

**Caution** : those **algebraic** K-groups will now be called $K_n(A)$. In order to avoid confusion, the topological K-groups for Banach algebras defined earlier will be denoted by $K_n^{top}(A)$.

As we have seen previously, homotopy groups play an important role in the definition and properties of the topological K-groups. Therefore, it is natural to look for an "algebraic" definition of them, especially for spaces like the general linear groups $GL_r(A)$ or their limit

$$GL(A) = \cup\, GL_r(A)$$

The first thing is to **define algebraically $\pi_0(GL(A))$**. In other words, one should say when two invertible matrices $\alpha_0$ and $\alpha_1$ are homotopic (or are in the same "algebraic" connected path component).

**DEFINITION**. *Two invertible matrices $\alpha_0$ and $\alpha_1$ are called **algebraically homotopic** if there exists an element $\alpha(x)$ of $GL(A[x])$ such that $\alpha(0) = \alpha_0$ and $\alpha(1) = \alpha_1$. This is an equivalence relation. The algebraically connected components form an abelian group denoted* by $K_1(A)$.

**Examples**. If F is a field or a Euclidean ring, it is easy to show (Gauss reduction) that $K_1(F) \cong F^*$, the group of invertible elements in F (take the kernel of the determinant map). One should note however that $K_1^{top}(\mathbf{C}) = 0$ and $K_1^{top}(\mathbf{R}) = \mathbf{Z}/2$.

A non trivial example [BMS] is the ring A of integers in a number field : one finds again $K_1(A) = A^*$. Note that $A^*$ has been determined by a theorem of Dirichlet (see [Sa] for a nice introduction to Number Theory).

For instance, if $A = \mathbf{Z}[x]/(1 + x + ... + x^{p-1})$, where p is an odd prime number, as in § 2, we have $K_1(A) \cong A^* \cong \mathbf{Z}^{(p-3)/2} \oplus \mathbf{Z}/p$.

An example of a ring A such that $K_1(A) \neq A^*$ is the following

$$A = \mathbf{Z}[x, y]/(x^2 + y^2 - 1)$$

One can show that the following matrix with coefficients in A :

$$\begin{pmatrix} x & -y \\ y & x \end{pmatrix}$$

is not algebraically homotopic to a 1 x 1 matrix, although its determinant is one (Hint : imbed A in the ring $C_\mathbf{R}(S^1)$).

Here is a table of $K_1$-groups which repeats the rings listed in the § 2

| A | $K_1(A)$ |
|---|---|
| **Z** | **Z**/2 |
| | The determinant of an invertible matrix with integral coefficients is ±1. Moreover, a matric of determinant 1 is a product of elementary matrices (use the Euclid algorithm). |
| $\mathbf{Z}(\sqrt{-5})$ | **Z**/2 |
| | This follows from a general theorem in Number Theory about the units in imaginary quadratic fields [Sa]. |
| $\mathbf{C}[G]$ | $(\mathbf{C}^*)^N$ |

Here G is a finite group, N the number of its irreducible representations. This computation follows from the fact that $\mathbf{C}[G]$ is a product of matrix algebras [Se].

| | |
|---|---|
| $C_\mathbf{C}(S^1)$ | $\mathrm{Map}(S^1, \mathbf{C}^*)$ |

In general, if A is a commutative Banach algebra, one has $K_1(A)$ isomorphic to $\pi_0(SL(A)) \oplus A^*$ (Proof in [M]).

| | |
|---|---|
| $C_\mathbf{R}(S^1)$ | $\mathrm{Map}(S^1, \mathbf{R}^*)$ |

Same comment.

| | |
|---|---|
| $C_\mathbf{C}(S^2)$ | $\mathrm{Map}(S^2, \mathbf{C}^*)$ |

Same comment.

| | |
|---|---|
| $C_\mathbf{R}(S^2)$ | $\mathrm{Map}(S^2, \mathbf{R}^*)$ |

Same comment.

| | |
|---|---|
| $M_n(A)$ | $K_1(A)$ |

This means that $K_1$-theory does not change when A is replaced by the ring of n x n matrices with coefficients in A (Morita invariance). Essentially the same proof as for the $K_0$-group.

| | |
|---|---|
| $\mathcal{K}(H)$ | 0 |

Essentially the same proof as for the $K_0$-group.

Coming back to the theory, it is easy to generalize this definition of $K_1(A)$ for rings without unit. If A is a k-algebra, following the previous definition for K(A), we define $K_1(A)$ as the kernel of the natural map

$$K_1(A^+) \longrightarrow K_1(k)$$

(it is easy to show that this definition is independent of k).

**Example.** An important example in Functional Analysis is the ring $\mathcal{K}$ of compact operators in an infinite dimensional Hilbert space. It is a challenging exercise to show that the $K_1$-group of this ring is reduced to 0.

We can push the analogy with topology further and introduce the "**loop ring**" $\Omega A$ of A as the non unital ring x(x-1) A[x], i.e. the two-sided ideal in A[x] consisting of polynomials P(x) such that P(0) = P(1) = 0. The topological analog of $\Omega A$ is $A(\mathbf{R}) = A(\,]0, 1[\,)$ in the previous section if A is a Banach algebra.

**DEFINITION**. *For* n ≥ 1, *the **algebraic** groups* $K_n(A)$ *are defined inductively by the following formula*

$$K_{n+1}(A) = K_n(\Omega A) = ... = K_1(\Omega^n A)$$

**Caution**. This definition does not coincide with Quillen's definition, except if the ring A is regular noetherian. This means that every finitely generated module admits a finite projective resolution. Many rings of geometric nature are of this type, for instance those associated to algebraic functions on affine algebraic varieties without singularities.

As for Banach algebras, a bilinear pairing

$$A \times C \longrightarrow B$$

induces a "**cup-product**"

$$K_n(A) \times K_p(C) \longrightarrow K_{n+p}(B)$$

An interesting example is A = C = B a commutative field (denoted by F). The above cup-product gives a "symbol map"

$$S : K_1(F) \otimes_{\mathbf{Z}} K_1(F) = F^* \otimes_{\mathbf{Z}} F^* \longrightarrow K_2(F)$$

**THEOREM** (Matsumoto). *The symbol map* S *is surjective. Its kernel is generated by the tensor products of the type* $x \otimes (1 - x)$, *where* $x \in F - \{0, 1\}$
[for a proof see [M1] ].

For example, in the case of the group $K_2(\mathbf{Q})$, where **Q** is the field of rational numbers, one finds

$$K_2(\mathbf{Q}) \cong \mathbf{Z}/2 \oplus [\bigoplus_p (\mathbf{Z}/p)^*]$$

where p runs through the set of all odd prime numbers.

Let us now state a striking conjecture [Ga] linking the groups $K_n(\mathbf{Z})$ with the **Bernoulli numbers** $B_k$. These rational numbers $B_k$ appear in the series expansion of the function $x/(1 - e^{-x})$. More precisely, one has

$$x/(1 - e^{-x}) = 1 + x/2 + \sum_{k=1}^{\infty} (-1)^{k-1} \frac{B_k}{(2k)!} x^{2k} + ...$$

Let us express the number $B_k$ as an irreducible fraction $c_k/d_k$. The conjecture is then the following :

$K_{8n}(\mathbf{Z}) = 0$ (for n > 0)
$K_{8n+1}(\mathbf{Z}) = \mathbf{Z} \oplus \mathbf{Z}/2$ (for n > 0)
$K_{8n+2}(\mathbf{Z}) = \mathbf{Z}/2 \oplus \mathbf{Z}/c_k$ (for k = 2n+1)
$K_{8n+3}(\mathbf{Z}) = \mathbf{Z}/8k.d_k$ (for k = 2n+1)
$K_{8n+4}(\mathbf{Z}) = 0$
$K_{8n+5}(\mathbf{Z}) = \mathbf{Z}$
$K_{8n+6}(\mathbf{Z}) = \mathbf{Z}/c_k$ (for k = 2n+2)
$K_{8n+7}(\mathbf{Z}) = \mathbf{Z}/4k.d_k$ (for k = 2n+2)

This conjecture has been verified for $K_i(\mathbf{Z})$ if i ≤ 6. Otherwise, we know at least that the right hand side is a lower bound for the groups $K_i(\mathbf{Z})$.

For instance, the conjecture predicts that $K_{23}(\mathbf{Z})$ is a cyclic group of order 65 520 (!). It is striking to see some partial periodicity appearing in the table as for **real** Banach algebras. Unfortunately, we are still far from understanding the full story. We recommend [Ga] for an overview and the relation with classical (unsolved) problems in Number Theory, like the Vandiver conjecture.

On the positive side, the first computations of algebraic K-groups were made by Quillen [Q1]. For a finite field with q elements (traditionnally called $\mathbf{F}_q$), he proved that

$K_{2i}(\mathbf{F}_q) = 0$ for i > 0 and
$K_{2i-1}(\mathbf{F}_q) \cong \mathbf{Z}/(q^i - 1)\mathbf{Z}$

From this result, we see that Bott periodicity cannot be true in a naïve sense. However, the situation is more favourable if one considers K-theory with finite coefficients $\mathbf{Z}/n$, say $K_i(A ; \mathbf{Z}/n)$. Such a theory has been introduced by W. Browder and the author [Br] [So] in the 70's : it fits into an exact sequence

$$K_i(A) \xrightarrow{.n} K_i(A) \longrightarrow K_i(A ; \mathbf{Z}/n) \longrightarrow K_{i-1}(A) \xrightarrow{.n} K_{i-1}(A)$$

If A is a field F, one has for instance $K_1(F ; \mathbf{Z}/n) \cong F^*/(F^*)^n$.

If F is an algebraically closed field, let us write $\mu_n$ for the group of $n^{th}$ roots of unity in F. It is isomorphic to $\mathbf{Z}/n$, but not canonically. Note that $K_*(F ; \mathbf{Z}/n)$ has a well defined ring structure in this case [Br] which is compatible with one shown before on $K_*(F)$.

**THEOREM** [Su]**.** *Let* F *be an algebraically closed field with a characteristic prime to* n. *Then*

$$K_{2i-1}(F ; \mathbf{Z}/n) = 0$$
$$K_{2i}(F ; \mathbf{Z}/n) \cong (\mu_n)^{\otimes i}$$

*Moreover, these isomorphisms are compatible with the ring structure. In other words, the cup-product*
$$K_{2i}(F ; \mathbf{Z}/n) \otimes K_2(F ; \mathbf{Z}/n) \longrightarrow K_{2i+2}(F ; \mathbf{Z}/n)$$
*is an isomorphism.*

This result is the true analog of Bott periodicity, since it implies that the K-groups with coefficients are periodic. Note however the existence of a twist since $\mu_n$ is not isomorphic to $\mathbf{Z}/n$ canonically.

**THEOREM** [Su]. *The obvious maps*
$$K_i(\mathbf{C} ; \mathbf{Z}/n) \longrightarrow K_i^{top}(\mathbf{C}; \mathbf{Z}/n)$$
$$K_i(\mathbf{R} ; \mathbf{Z}/n) \longrightarrow K_i^{top}(\mathbf{R}; \mathbf{Z}/n)$$

*are isomorphisms. In particular, the former groups are periodic of period 2 and the latter are periodic of period 8.*

At this point, one may ask how to compute $K_i(F ; \mathbf{Z}/n)$ for an arbitrary field (the groups $K_*(F)$ themselves seem untractable for the moment). The general feeling among specialists (not yet completely understood) is that we should calculate them from the groups $K_*(\overline{F} ; \mathbf{Z}/n)$, where $\overline{F}$ is the separable closure of F, taking into account the Galois action on them : this requires a lot of machinery from Algebra and Algebraic Topology.

**Applications.** One can find many applications of Algebraic K-theory in the references [R][Sr]. We just mention two of them.

The first one is relevant to Differential Topology. Let M be a sphere $S^n$. The space $\mathbb{P}(M)$ of pseudo-isotopies of M is the space of diffeomorphisms of the product $M \times [0, 1]$ which induce the identity on $M \times \{0\}$. By a result of Waldhausen [L], the homotopy groups $\pi_i(\mathbb{P}(M)) \otimes_{\mathbf{Z}} \mathbf{Q}$ are isomorphic to $K_{i+2}(\mathbf{Z}) \otimes_{\mathbf{Z}} \mathbf{Q}$. In particular, according to Borel [Bo], these groups are periodic of period 4, isomorphic to 0, 0, $\mathbf{Q}$, 0, 0, 0, $\mathbf{Q}$... (starting from i = 0).

Another spectacular application is due to Merkujev and Suslin [MS] and involves a well known algebraic invariant called the Brauer group of a field F. One considers central simple algebras A over F. Two such algebras A and A' are called equivalent if $M_p(A)$ is isomorphic to $M_p(A')$ for well choosen p and p'. The Brauer group is the quotient of this set of algebras by the equivalence relation (it is a group by the law induced by the tensor product of algebras). If the field F contains the $n^{th}$ roots of unity, then the n-torsion of the Brauer group is isomorphic to $K_2(F)/nK_2(F)$. For instance, if n = 2, this result implies that any element of the 2-torsion of the Brauer group can be represented by the tensor product of finitely many quaternion algebras.

# SELECTED REFERENCES

## Books on K-theory and applications

## Books on Algebraic Topology

## More books and papers

**Max Karoubi, Université Paris 7**
**2, Place Jussieu, 75251 Paris (France)**